\documentclass[11pt,oneside]{article}
 \usepackage{afterpage}
 \usepackage{fancyhdr}
 \usepackage{lipsum}
 \newcommand\shorttitle{Flag curvature of a homogeneous generalized Kropina metric}
 \newcommand\authors{g. shanker, j. kaur, seema}
 
 \usepackage{supertabular, setspace,hyperref}
 \usepackage[a4paper, total={6.5in, 9.5in}]{geometry}
 \usepackage[utf8]{inputenc}
 \usepackage{amsmath}
 \usepackage{mathtools}
 \usepackage{amsfonts}
 \usepackage{amssymb}
 \usepackage{amsthm}
 \usepackage{caption}
 \usepackage{subcaption}
 \usepackage{authblk}
 \usepackage[sort]{cite}
 \newtheorem{theorem}{Theorem}[section]

 \newtheorem{example}[theorem]{Example}

 \newtheorem{remark}{\sc Remark}
 \newtheorem{lemma}{\sc Lemma}[section]
 \newtheorem{corollary}{\sc Corollary}[section]
 \newtheorem{definition}{\sc Definition}[section]

 \newcommand{\be}{\begin{eqnarray}}
 	\newcommand{\ee}{\end{eqnarray}}
 \newcommand{\Be}{\begin{eqnarray*}}
 	\newcommand{\Ee}{\end{eqnarray*}}
 \newcommand{\bee}{\begin{equation}}
 	\newcommand{\eee}{\end{equation}}
 \newcommand{\ba}{\begin{array}}
 	\newcommand{\ea}{\end{array}}
 \newcommand{\bl}{\begin{lemma}}
 	\newcommand{\el}{\end{lemma}}
 \newcommand{\bd}{\begin{definition}}
 	\newcommand{\ed}{\end{definition}}
 \newcommand{\bt}{\begin{theorem}}
 	\newcommand{\et}{\end{theorem}}
 \newcommand{\bp}{\begin{proof}}
 	\newcommand{\ep}{\end{proof}}
 \newcommand{\bi}{\begin{itemize}}
 	\newcommand{\ei}{\end{itemize}}
 \newcommand{\br}{\begin{remark}}
 	\newcommand{\er}{\end{remark}}
 \newcommand{\bc}{\begin{corollary}}
 	\newcommand{\ec}{\end{corollary}}
 \newcommand{\bex}{\begin{example}}
 	\newcommand{\eex}{\end{example}}

 \usepackage{chngcntr}
 
 \counterwithin*{equation}{section}
 \counterwithin*{equation}{section}
 
\begin{document}
 	\afterpage{\cfoot{\thepage}}
 	\clearpage
 	\date{}

 	 \title	{\textbf{On the flag curvature of a homogeneous Finsler space with generalized $m$-Kropina metric}}
 	\author[]{Gauree Shanker}
 	\author[]{Jaspreet Kaur\thanks{corresponding author, Email: gjaspreet303@gmail.com}}
 	\author[]{Seema} 
 	\affil[]{\footnotesize Department of Mathematics and Statistics,
 		Central University of Punjab, Bathinda, Punjab-151401, India\\
 		Email:  grshnkr2007@gmail.com, gjaspreet303@gmail.com, seemajangir2@gmail.com}
 	\maketitle
 	\begin{center}
 		\textbf{Abstract}
 	\end{center}
 	\begin{small}
 		In this paper, first, we give an explicit formula for the flag curvature of a homogeneous Finsler space with generalized $m$-Kropina metric. Then, we show that, under a mild condition, the
 		two definitions of naturally reductive homogeneous Finsler space are
 		equivalent for aforesaid metric. Finally, we study the flag curvature of naturally reductive homogeneous Finsler spaces with generalized $m$-Kropina metric.
 	\end{small}
 	\section{Introduction}
Finsler geometry is one of the important reasearch area in differential geometry which has been
developed very rapidly in recent years. One reason for this development is its application in many areas of natural science such as biology and physics \cite{R2,R3}. Finsler geometry is just Riemannian geometry without quadratic restriction \cite{R4}. In $1972$, M. Matsumoto \cite{R9} introduced the notion of $(\alpha,\beta)$-metrics in Finsler geometry. Some important $(\alpha,\beta)$-metrics are Randers metric, Matsumoto metric, Kropina metric, generalized $m$-Kropina metric, square metric, etc. Many authors \cite{R12,R56,R10,R8,R18,R23} have studied various properties of $(\alpha,\beta)$-metrics. Generalized $m$-Kropina metric belongs to the large class of $(\alpha,\beta)$-metrics. Kropina metric and generalized $m$-Kropina metric, both have a lot of applications in other branches of science such as physics, irreversible thermodynamics, electron optics with a magnetic field, etc \cite{R27,R26}.\\
$(\alpha,\beta)$-metric on a connected smooth $n$-manifold $M$ is a Finsler metric which can be written in the form $$F=\alpha\phi(s);~s=\dfrac{\beta}{\alpha},$$ \\
 where $\alpha:=\alpha(x,y)=\sqrt{a_{ij}(x)y^{i}y^{j}}$ is a Riemannian metric and $\beta:=\beta(y)=b_{i}(x)y^{i}$ is a $1$-form on $M$.\\
  In particular, if $\phi(s)=\dfrac{1}{s^{m}}(m\neq0,-1)$, then the Finsler metric $$F=\dfrac{\alpha^{m+1}}{\beta^{m}}$$ is called generalized $m$-Kropina metric.
\\
Flag curvature is the most important quantity in Finsler geometry because it is a generalization of the sectional curvature of Riemannian metric. In general, the computation of the flag curvature of a Finsler metric is difficult, therefore it is very important to find an explicit and applicable formula for the flag curvature. Many authors \cite{R14, R15,R17,R20,R22} have worked in this area.\\
The notion of naturally reductive Riemannian metric was first introduced by Kobayashi and Namizu \cite{R21}. The naturally reductive spaces have been investigated by several authors as a natural generalization of Riemannian symmetric spaces. In literature, there are two versions of the definition of naturally reductive spaces. The first definition was
given by Deng and  Hou \cite{R25}. In this definition, authors have supposed that the metric should be Berwaldian. The second one was given by Latifi \cite{R24}. Deng and Hou, in \cite{R10} have proved that if a homogeneous Finsler space is naturally reductive by Latifi's definition, then it must be naturally reductive by Deng and Hou's definition and Berwaldian. Parhizkar and Moghaddam \cite{R28} prove that both the
definitions of naturally reductive homogeneous Finsler spaces are equivalent under
the consideration of a mild condition.\\
This paper is organized as follows:\\
 Section $2$ includes basic information on homogeneous Finsler spaces needed in this paper. In section $3$, we obtain the formula for flag curvature of a homogeneous Finsler space with generalized Kropina metric under certain conditions. In section $4$, we show that, under a mild condition, two definitions of naturally reductive homogeneous Finsler space are equivalent for aforesaid metric. Finally, in section $5$, we derive an explicit formula for flag curvature of a naturally reductive homogeneous Finsler space with generalized $m$-Kropina metric.
	\section{Preliminaries}
   In this section, we give some basic definitions and results of Finsler spaces and naturally reductive homogeneous spaces. We refer \cite{R1,R5,R6} for notations and symbols.
    	\begin{definition}
    		Let $V$ be a real vector space of dimension $n$. A smooth function $F:V\rightarrow [0,\infty)$ is called a Minkowski norm if 
    	\begin{itemize}
    		\item[(i)] $F(\lambda y)=\lambda F(y) ~\forall~\lambda>0,~y\in V\setminus\{0\},$ i.e., $F$ is positively homogeneous.
    	\item[(ii)] For each fixed $y\in V\setminus\{0\},$ the bilinear function
    	 $g_{_{Y}}:V\times V\rightarrow\mathbb{R}$ defined by 
    $$g_{_{Y}}(u,v)=\dfrac{1}{2}\dfrac{\partial^{2}}{\partial s\partial t}F^{2}(y+su+tv)\Big|_{s=t=0}$$ is positive-definite at every point of $V\setminus\{0\}$.
    	\end{itemize}
     Minkoswki space is a linear vector space $V$ equipped with a Minkowski norm $F$, denoted by $(V,F).$  \\
     For any basis $\{u_{1},u_{2},...,u_{n}\}$ of $V$ and $y=y^{i}u_{i}\in V$, the Minkowski norm $F$ can be written as $F(y)=F(y^{1},y^{2},...,y^{n})$ and the Hessian matrix is $(g_{ij})=\bigg[\dfrac{1}{2}F^{2}\bigg]_{y^{i}y^{j}}$.
    	\end{definition}
   
    \begin{definition}
    	A Finsler metric on a smooth manifold $M$ is a real valued function $F:TM\rightarrow[0,\infty)$ such that
    	\begin{itemize}
    		\item[(i)] $F$ is smooth on the slit tangent bundle $TM\setminus\{0\}$.
    		\item[(ii)] The restriction of $F$ to each tangent space $T_{x}M,~x\in M$, is a Minkowski norm.
    	\end{itemize} 
    \end{definition}
 \begin{definition}
	A Finsler space $(M,F)$ is called a Berwald space if Chern connection of $(M,F)$ is a linear connection on $TM$. Equivalently, if each of the Chern connection coefficient $\Gamma^{i}_{jk}$, in natural standard coordinate system, have no $y$ dependence, then Finsler space $(M,F)$ becomes a Berwald space.
\end{definition}
\begin{definition}
	Let $(M,F)$ be a Finsler space. The flag curvature $K=K(P,Y)$ is a function of tangent planes $P=span\{Y,U\}\subset T_{x}M$ and directions $Y\in T_{x}M\setminus\{0\}$. The pair $(P,Y)$ is called a flag and $Y$ is called the flag pole. The flag curvature is defined by 
	\begin{equation}\label{p}
	K(P,Y)=\dfrac{g_{_{Y}}(U,R(U,Y)Y)}{g_{_{Y}}(Y,Y)g_{_{Y}}(U,U)-g_{_{Y}}^{2}(Y,U)}.
	\end{equation}
	 \quad If $F$ is Riemannian, $K=K(P)$ is independent of $Y\in P\setminus\{0\}$ and is called the sectional curvature in Riemannian geometry. If $K=K(x,Y)$ is a scalar function on the slit tangent bundle $TM\setminus\{0\}$. then Finsler metric $F$ is said to be of scalar curvature.
\end{definition}

\begin{lemma}
	Let $F=\alpha\phi(s);~s=\dfrac{\beta}{\alpha}$, where $\alpha$ is a Riemannian metric, $\beta$ is a $1$-form with $||\beta||_{\alpha}<b_{0}$ and $\phi$ is a smooth function on an open interval $(-b_{0},b_{0})$. Then $F$ is Finsler metric if and only
	if $\phi$ satisfies the following condition :
	$$\phi(s)>0, ~\phi(s)-s\phi'(s)+(b^{2}-s^{2})\phi''(s)>0~\forall~|s|\leq b< b_{0}.$$
\end{lemma}
\begin{definition}
	Let $G$ be a smooth manifold with structure of an abstract group. If the map $\varphi: G\times G\rightarrow G,$ defined by $\varphi(g_{1},g_{2})=g_{1}g_{2}^{-1}$ is smooth, then $G$ is called a Lie group.
\end{definition} 
	If $G$ is a Lie group and $M$, a smooth manifold, then smooth action of $G$ on $M$ is a smooth map $\varphi:G\times M\rightarrow M$ satisfying $\varphi(g_{2},\varphi(g_{1},x))=\varphi(g_{2}g_{1},x)~~\forall~ g_{1},g_{2}\in G,~x\in M$.\\
     If a Lie group acts smoothly on a smooth manifold $M$, then $G$ is called a Lie transformation group of $M$.
\begin{definition}
	A Finsler space $(M,F)$ is called homogeneous if the group $I(M,F)$ of isometries acts transitively on $M$. It can be naturally identified with the
	homogeneous space $G/H$, where $H$ is the isotropy group at $p\in M$.
\end{definition}
\begin{definition}\cite{R21}
Let $G/H$ be a homogeneous space with an invariant Riemannian metric $\tilde{g}$ and $\mathfrak{g}$ and $\mathfrak{h}$ are the Lie algebras of $G$ and $H$ respectively. Then $G/H$ is called naturally reductive if there exists an $Ad(H)$-invariant decomposition $\mathfrak{g}=\mathfrak{h}+\mathfrak{m}$ such that
$$\langle[X,Y]_{\mathfrak{m}},Z\rangle+\langle X,[Y,Z]_{\mathfrak{m}}\rangle=0 ~\forall~ X,Y,Z\in\mathfrak{m},$$
where $\mathfrak{m}$ is subspace of $\mathfrak{g}$ such that $Ad(h)\mathfrak{m}\subset\mathfrak{m}~ \forall h\in H$  and $\langle~,~\rangle$ is bilinear form on $\mathfrak{m}$ induced by $\tilde{g}$ and $[~,~]_{\mathfrak{m}}$ is projection of $[~,~]$ to $\mathfrak{m}$.
\end{definition}
\begin{definition}\cite{R25}
	A homogeneous Finsler space $G/H$ with an invariant Finsler metric
	$F$ is said to be naturally reductive if there exists an invariant Riemannian metric $\tilde{g}$ on $G/H$ such that $(G/H, \tilde{g})$ is naturally reductive and the Chern connection of $F$
	coincides with the Levi-Civita connection of $\tilde{g}$.
\end{definition}
\begin{definition}\cite{R24}
	A homogeneous space $G/H$ with an invariant Finsler metric is called naturally reductive if there exists an $Ad(H)$-invariant decomposition $\mathfrak{g}=\mathfrak{h}+\mathfrak{m}$ such that $$g_{_{Y}}([x,u]_{\mathfrak{m}},v)+g_{_{Y}}(x,[u,v]_{\mathfrak{m}})+2C_{y}([x,y]_{\mathfrak{m}},u,v)=0,$$ where $y(\neq0),x,u,v\in\mathfrak{m}$.
\end{definition}
\begin{lemma}\label{ER}\cite{R24}
	Let $G/H$ be a naturally reductive homogeneous space with an $Ad(H)$-invariant decomposition $\mathfrak{g}=\mathfrak{h}+\mathfrak{m}$ and a $G$-invariant indefinite Riemannian metric $\tilde{g}$. Then the curvature tensor $R$ of the Riemannian connection satisfies 
	\begin{equation*}
	R(X,Y)Y=\dfrac{1}{4}\big[Y,[X,Y]_{\mathfrak{m}}\big]_{\mathfrak{m}}+\big[Y,[X,Y]_{\mathfrak{h}}\big] ~\forall~X,Y\in\mathfrak{m}.
	\end{equation*}
	\end{lemma}

\section{Flag curvature of homogeneous Finsler spaces}
In this section, we derive an explicit formula for the flag curvature of a homogeneous Finsler space with generalized $m$-Kropina metric. For this purpose, we use  P$\ddot{u}$ttann's \cite{R7} formula for curvature tensor of invariant metric $\langle~,~\rangle$ on compact homogeneous space.\\
Let $G$ be a compact Lie group, $H$ be a closed subgroup of $G$ with Lie algebras $\mathfrak{g}$ and $\mathfrak{h}$ respectively and $\tilde{g}$ be a bi-invariant Riemannian metric on $G$. The tangent space of the homogeneous space is given by orthogonal complement $\mathfrak{m}$ of $\mathfrak{h}$ in $\mathfrak{g}$. Each invariant metric $\tilde{g}$ on $G/H$ is determined by its restriction to $\mathfrak{m}$. Ad$(H)$-invariant inner product on $\mathfrak{m}$ can be extended to an Ad$(H)$-invariant inner product on $\mathfrak{g}$ by taking $\tilde{g_{0}}$ for the components of $\mathfrak{m}$. In this way, the metric on $G/H$ determines a unique left invariant metric $\tilde{g}$ on $G$. Also, at the identity $e$ of $G$, the value of $\tilde{g_{0}}$ and $\tilde{g}$ are inner products on $g$, we denote them by $\langle\langle~,~\rangle\rangle$ and $\langle~,~\rangle$ respectively. Moreover, $\langle~,~\rangle$ determines a positive definite endomorphism $\varphi$ of $\mathfrak{g}$ such that $\langle X,Z\rangle=\langle\langle\varphi(X),Z\rangle\rangle~\forall~X,Z\in\mathfrak{g}$.\\
 P$\ddot{u}$ttann's formula for curvature tensor of invariant metric $\langle~,~\rangle$ on compact homogeneous space $G/H$ is given by 
\begin{equation*}
\begin{split}
\langle R(X,Y)Z,W\rangle&=\dfrac{1}{2}\Big\{\langle\langle B_{-}(X,Y),[Z,W]\rangle\rangle+\langle\langle [X,Y],B_{-}(Z,W)\rangle\rangle\Big\}\\
&+\dfrac{1}{4}\Big\{\langle[X,W],[Y,Z]_{\mathfrak{m}}\rangle-\langle[X,Z],[Y,W]_{\mathfrak{m}}\rangle-2\langle[X,Y],[Z,W]_{\mathfrak{m}}\rangle\Big\}\\
&+\langle\langle B_{+}(X,W),\varphi^{-1}B_{+}(Y,Z)\rangle\rangle-\langle\langle B_{+}(X,Z),\varphi^{-1}B_{+}(Y,W)\rangle\rangle,
\end{split}
\end{equation*}
where
$$B_{+}(X,Y)=\dfrac{1}{2}\Big([X,\varphi Y]+[Y,\varphi X]\Big),$$
and
$$B_{-}(X,Y)=\dfrac{1}{2}\Big([\varphi X,Y]+[X,\varphi Y]\Big)$$ are the bilinear symmetric and skew-symmetric maps respectively.

 \begin{theorem}\label{1}
    	Let $H$ be closed subgroup of compact Lie group $G$ with $\mathfrak{g}$ and $\mathfrak{h}$, Lie algebras of $G$ and $H$ respectively. Further, let $\tilde{g_{0}}$ be a bi-invariant metric on Lie group $G$, $\tilde{g}$ be invariant Riemannian metric on homogeneous space $G/H$  and  $\tilde{X}$, an invariant vector field on $G/H,$ which is parallel with respect to $\tilde{g}$ such that $\sqrt{\tilde{g}(\tilde{X},\tilde{X})}<1$ and $\tilde{X}_{H}=X$. Let $F=\dfrac{\alpha^{m+1}}{\beta^{m}}~(m\neq 0,-1)$ be a generalized $m$-Kropina metric arising by $\tilde{g}$ and $\tilde{X}.$ Let $(P,Y)$ be a flag in $T_{H}(G/H)$ such that $\{U,Y\}$ is an orthonormal basis of $P$ with respect to $\langle~,~\rangle.$ Then the flag curvature of the flag $(P,Y)$ in $T_{H}(G/H)$ is given by 
    	\begin{equation}\label{i}
        K(P,Y)=\dfrac{\langle X, Y\rangle^{2m}\big[(m+1)\langle X,Y\rangle^{2}\langle U,R(U,Y)Y\rangle+m(2m+1)\langle X,U\rangle \langle X,R(U,Y)Y\rangle\big]}{(m+1)\big[m\langle X,U\rangle^{2}+\langle X,Y\rangle^{2}\big]},
    	 \end{equation}
    	 where 
    	 \begin{equation}\label{j}
    	 \begin{split}
    	  \langle X,R(U,Y)Y\rangle&=-\dfrac{1}{4}\bigg(\langle[\varphi U,Y]+[U,\varphi Y],[Y,X]\rangle+\langle[U,Y],[\varphi Y,X]+[Y,\varphi X]\rangle\bigg) \\
    	  &-\dfrac{3}{4}\langle[Y,U],[Y,X]_{\mathfrak{m}}\rangle-\dfrac{1}{2}\langle[U,\varphi X]+[X,\varphi U],\varphi^{-1}([Y,\varphi Y])\rangle\\
    	  &+\dfrac{1}{4}\langle[U,\varphi Y]+[Y,\varphi U],\varphi^{-1}([Y,\varphi X]+[X,\varphi Y])\rangle,
    	 \end{split}  	
    	 \end{equation}
    	 and 
    	 \begin{equation}\label{k}
    	 \begin{split}
    	  \langle U,R(U,Y)Y\rangle&=\dfrac{1}{2}\big(\langle[\varphi U,Y]+[U,\pi Y],[Y,U]\rangle\big)\\ &+\dfrac{3}{4}\langle[Y,U],[Y,U]_{m}\rangle+\langle[U,\varphi U],\varphi^{-1}([Y,\varphi Y])\rangle \\
    	 &-\dfrac{1}{4}\langle[U,\varphi Y]+[Y,\varphi U],\varphi^{-1}([Y,\varphi U]+[U,\varphi Y])\rangle.
    	 \end{split}    	 
    	 \end{equation}
    	 \textit{Proof.} Since $\tilde{X}$ is parallel with respect to $\tilde{g}$, $\beta$ is parallel with respect to $\tilde{g}$. Therefore, the Chern connection of $F$ coincide with the Levi-Civita connection of $\tilde{g}$. Thus the Finsler metric $F$ and the Riemannian metric $\tilde{g}$ have same curvature tensor and we denote it by $R$. Therefore,
    	 we have 
    	 \begin{equation*}
    	 F(Y)=\dfrac{(\sqrt{\langle Y,Y\rangle})^{m+1}}{\langle X,Y\rangle^{m}}.
    	 \end{equation*}
     By using the definition of $g_{_{Y}}(U,V)$, after some computations, we get
    	 \begin{equation}\label{l}
    	 \begin{split}
    	  g_{_{Y}}(U,V)&=\dfrac{\langle Y,Y\rangle^{m-1}}{\langle X,Y\rangle^{2m+2}}\Big[2m(m+1)\langle X,Y\rangle^{2}\langle U,Y\rangle\langle Y,V\rangle-2m(m+1)\langle X,Y\rangle\langle Y,Y\rangle\langle Y,U\rangle\langle X,V\rangle\\
    	  &-2m(m+1)\langle X,Y\rangle\langle Y,Y\rangle\langle X,U\rangle\langle Y,V\rangle+(m+1)\langle X,Y\rangle^{2}\langle Y,Y\rangle\langle U,V\rangle\\
    	  &+m(2m+1)\langle Y,Y\rangle^{2}\langle X,U\rangle\langle X,V\rangle\Big].
    	 \end{split}
    	   \end{equation}
    	   As $\{U,Y\}$is an orthonormal basis for $P$ with respect to $\langle~,~\rangle$, the equation \eqref{l} can be written as:
    	   \begin{equation}\label{m}
    	   \begin{split}
    	    g_{_{Y}}(U,V)&=\dfrac{1}{\langle X,Y\rangle^{2m+2}}\Big[(m+1)\langle X,Y\rangle^{2}\langle U,V\rangle-2m(m+1)\langle X,Y\rangle\langle X,U\rangle\langle Y,V\rangle\\
    	    &+m(2m+1)\langle X,U\rangle\langle X,V\rangle\Big].
    	   \end{split}
    	   \end{equation}
    	  Therefore,
    	    \begin{equation}\label{o}
    	   \begin{split}
    	   g_{_{Y}}(U,R(U,Y)Y)&=\dfrac{1}{\langle X,Y\rangle^{2m+2}}\Big[m\langle X,Y\rangle^{2}\langle U,R(U,Y)Y\rangle-2m(m+1)\langle X,Y\rangle\langle X,U\rangle\langle Y,R(U,Y)Y\rangle\\
    	   &+m(2m+1)\langle X,U\rangle\langle X,R(U,Y)Y\rangle\Big].
    	   \end{split}
    	   \end{equation}
    	    \vspace{0.3cm}
    	   From equation \eqref{m}, we get following three equations:\\
    	   $$g_{_{Y}}(Y,Y)=\dfrac{1}{\langle X,Y\rangle^{2m}},$$
    	   $$g_{_{Y}}(U,Y)=- \dfrac{m\langle X,U\rangle}{\langle X,Y\rangle^{2m+1}},$$
    	  and
    	   $$g_{_{Y}}(U,U)=\dfrac{1}{\langle X,Y\rangle^{2m+2}}\Big[(m+1) \langle X,Y\rangle^{2}+m(2m+1)\langle X,U\rangle^{2}\Big].$$
    	   From above three equations, we get
    	   \begin{equation}\label{n}
    	   g_{_{Y}}(Y,Y)g_{_{Y}}(U,U)-g_{_{Y}}^{2}(U,Y)=\dfrac{(m+1)}{\langle X,Y\rangle^{4m+2}}\Big[m\langle X,U\rangle^{2}+\langle X,Y\rangle^{2}\Big].
    	   \end{equation}
    	  Finally, using equations \eqref{j}, \eqref{k}, \eqref{o} and \eqref{n}  in equation \eqref{p}, we get the required proof.
    \end{theorem}
    \section{Naturally reductive homogeneous Finsler space}
    In this section, we prove that if a homogeneous Finsler space with generalized $m$-Kropina metric is naturally reductive in the sense of Latifi, under a mild condition, then it is naturally reductive in the sense of Deng and Hou and vice-versa.
   
    \begin{theorem}
	Let $G$ be a compact Lie group and $H$ be its closed subgroup with Lie algebras $\mathfrak{g}$ and $\mathfrak{h}$ respectively. Also, let $(G/H,F)$ be a homogeneous Finsler space of Berwald type with an invariant Riemannian metric $\langle ~,~\rangle$ and an invariant vector field $X$ such that $\tilde{X}(H)=X$. Then $(G/H,F)$ is naturally reductive if and only if the Riemannian space $(G/H,\langle~,~\rangle)$ is naturally reductive.\\
    \textit{Proof.}	Let $Y(\neq 0),~Z\in\mathfrak{m}$. From equation \eqref{l}, we have
    \begin{equation*}
    \begin{split}
    g_{_{Y}}(Y,[Y,Z]_{\mathfrak{m}})&=\dfrac{\langle Y,Y\rangle^{m-1}}{\langle X,Y\rangle^{2m+2}}\Bigg[2m(m+1)\langle X,Y\rangle^{2}\langle Y,Y\rangle\langle Y,[Y,Z]_{\mathfrak{m}}\rangle\\
    &-2m(m+1)\langle X,Y\rangle\langle Y,Y\rangle^{2}\langle X,[Y,Z]_{\mathfrak{m}}\rangle
    -2m(m+1)\langle X,Y\rangle^{2}\langle Y,Y\rangle\langle Y,[Y,Z]_{\mathfrak{m}}\rangle\\
    &+(m+1)\langle X,Y\rangle^{2}\langle Y,Y\rangle\langle Y,[Y,Z]_{\mathfrak{m}}\rangle+m(2m+1)\langle Y,Y\rangle^{2}\langle X,Y\rangle\langle X,[Y,Z]_{\mathfrak{m}}\rangle\Bigg],
   \end{split}
    \end{equation*}
    i.e.,
    \begin{equation}\label{q}
    g_{_{Y}}(Y,[Y,Z]_{\mathfrak{m}}) =\dfrac{\langle Y,Y\rangle^{m}}{\langle X,Y\rangle^{2m+1}}\Bigg[(m+1)\langle X,Y\rangle\langle Y,[Y,Z]_{\mathfrak{m}}\rangle-m\langle Y,Y\rangle\langle X,[Y,Z]_{\mathfrak{m}}\rangle\Bigg].
    \end{equation}
    Since Chern connection of $(G/H,F)$ coincide with the Levi-Civita connection of $(G/H,\langle~,~\rangle)$.\\
     From equation \eqref{q}, we get 
     \begin{equation}\label{s}
     \langle X,[Y,Z]_{\mathfrak{m}}\rangle=0 ~ \forall~Z\in\mathfrak{m}.
     \end{equation}
     Now, suppose $(G/H,F)$ is naturally reductive, then we have
     \begin{equation*}
     g_{_{Y}}([Z,U]_{\mathfrak{m}},V)+g_{_{Y}}([Z,V]_{\mathfrak{m}},U)+2 C_{Y}([Z,Y]_{\mathfrak{m}},U,V)=0~ \forall~Y\neq 0,Z,U,V\in\mathfrak{m}.
     \end{equation*}
     Therefore, we can write
     \begin{equation*}
      g_{_{Y}}([Y,U]_{\mathfrak{m}},V)+g_{_{Y}}([Y,V]_{\mathfrak{m}},U)+2 C_{Y}([Y,Y]_{\mathfrak{m}},U,V)=0 ~\forall~Y\neq 0,
     \end{equation*}
     i.e.,
     \begin{equation}\label{t}
     g_{_{Y}}([Y,U]_{\mathfrak{m}},V)+g_{_{Y}}([Y,V]_{\mathfrak{m}},U)=0,
     \end{equation}
     i.e.,
     \begin{equation}\label{u}
      g_{_{Y}}([Y,Z]_{\mathfrak{m}},Y)=0.
     \end{equation}
     Using equations \eqref{s}, \eqref{u} in \eqref{q}, we get
     \begin{equation}\label{v}
     \langle[Y,Z]_{\mathfrak{m}},Y\rangle=0.
     \end{equation}
     Also, from the equation \eqref{l}, we can write
     \begin{equation}\label{w}
    \begin{split}
    g_{_{Y}}([Y,U]_{\mathfrak{m}},V)=\dfrac{\langle Y,Y\rangle^{m-1}}{\langle X,Y\rangle^{2m+2}}&\Bigg[2m(m+1)\langle X,Y\rangle^{2}\langle Y,V\rangle\langle [Y,U]_{\mathfrak{m}},Y\rangle\\
    &-2m(m+1)\langle X,Y\rangle\langle Y,Y\rangle\langle Y,[Y,U]_{\mathfrak{m}}\rangle\langle X,V\rangle\\
    &-2m(m+1)\langle X,Y\rangle\langle Y,Y\rangle\langle X,[Y,U]_{\mathfrak{m}}\rangle\langle Y,V\rangle\\
    &+(m+1)\langle X,Y\rangle^{2}\langle Y,Y\rangle\langle[Y,U]_{\mathfrak{m}},V\rangle\\
   &+m(2m+1)\langle Y,Y\rangle^{2}\langle X,V\rangle\langle X,[Y,U]_{\mathfrak{m}}\rangle\Bigg].\\
    \end{split}
     \end{equation}
     Using equations \eqref{w}, \eqref{s} and \eqref{u}, we get
     \begin{equation}\label{x}
      g_{_{Y}}([Y,U]_{\mathfrak{m}},V)=\dfrac{(m+1)\langle Y,Y\rangle^{m}}{\langle X,Y\rangle^{2m}}\langle[Y,U]_{\mathfrak{m}},V\rangle.
     \end{equation}
     Similary, we get 
     \begin{equation}\label{y}
     g_{_{Y}}([Y,V]_{\mathfrak{m}},U)=\dfrac{(m+1)\langle Y,Y\rangle^{m}}{\langle X,Y\rangle^{2m}}\langle[Y,V]_{\mathfrak{m}},U\rangle.
     \end{equation}
     From equations \eqref{t}, \eqref{x} and \eqref{y}, we have
     \begin{equation*}
     \langle[Y,U]_{\mathfrak{m}},V\rangle+\langle[Y,V]_{\mathfrak{m}},U\rangle=0.
     \end{equation*}
     Hence $(G/H,\langle~,~\rangle)$ is naturally reductive.\\
     Conversely, let $(G/H,\langle~,~\rangle)$ be naturally reductive.
     From equations \eqref{l} and \eqref{s}, we can write
  
     \begin{equation}\label{ju}
     \begin{split}
      g_{_{Y}}(\langle[Z,U]_{\mathfrak{m}},V)&=\dfrac{\langle Y,Y\rangle^{m-1}}{\langle X,Y\rangle^{2m+1}}\Bigg[2m(m+1)\langle[Z,U]_{\mathfrak{m}},Y\rangle\bigg(\langle X,Y\rangle\langle Y,V\rangle-2m(m+1)\langle X,V\rangle\langle Y,Y\rangle\bigg)\\
      &+(m+1)\langle X,Y\rangle\langle Y,Y\rangle\langle [Z,U]_{\mathfrak{m}},V\rangle\Bigg].
     \end{split}
     \end{equation} 
     Similary, we have 
     \begin{equation}\label{jt}
      \begin{split}
     g_{_{Y}}(\langle[Z,V]_{\mathfrak{m}},U)&=\dfrac{\langle Y,Y\rangle^{m-1}}{\langle X,Y\rangle^{2m+1}}\Bigg[\langle[Z,V]_{\mathfrak{m}},Y\rangle\bigg(2m(m+1)\langle X,Y\rangle\langle Y,U\rangle-\langle X,U\rangle\langle Y,Y\rangle\bigg)\\
     &+(m+1)\langle X,Y\rangle\langle Y,Y\rangle\langle [Z,V]_{\mathfrak{m}},U\rangle\Bigg].
     \end{split}
     \end{equation}
     Also, Cartan tensor is given by 
     \begin{equation*}
     C_{Y}(Z,U,V)=\dfrac{1}{2}\dfrac{d}{dt}\Big[g_{_Y+tV}(Z,U)\Big]\bigg|_{t=0}.
     \end{equation*}
     After some calculation, we get
     \begin{equation}{\label{jh}}
     \begin{split}
      2C_{Y}([Z,Y]_{\mathfrak{m}},U,V)&=\dfrac{2m(m+1)\langle Y,Y\rangle^{m-1}}{\langle X,Y\rangle^{2m+2}}\Bigg[\langle[Z,Y]_{\mathfrak{m}},V\rangle\bigg(\langle X,Y\rangle^{2}\langle Y,U\rangle-\langle X,Y\rangle\langle X,U\rangle\langle Y,Y\rangle\bigg)\\
     &+\langle[Z,Y]_{\mathfrak{m}},U\rangle\bigg(\langle X,Y\rangle^{2}\langle Y,V\rangle-\langle X,Y\rangle\langle X,V\rangle\langle Y,Y\rangle\bigg)\Bigg].
     \end{split}    
     \end{equation}
     Adding equations \eqref{ju}, \eqref{jt} and \eqref{jh}, we get 
     \begin{equation*}
     \begin{split}
      &g_{_{Y}}(\langle[Z,U]_{\mathfrak{m}},V)+ g_{_{Y}}(\langle[Z,V]_{\mathfrak{m}},U)+ 2C_{Y}([Z,Y]_{\mathfrak{m}},U,V)\\
      &=\dfrac{\langle Y,Y\rangle^{m-1}}{\langle X,Y\rangle^{2m+2}}\Bigg[(m+1)\langle X,Y\rangle^{2}\langle Y,Y\rangle
      \bigg(\langle[Z,V]_{\mathfrak{m}},U\rangle+\langle[Z,U]_{\mathfrak{m}},V\rangle\bigg)\\
      &+2m(m+1)\langle X,Y\rangle\bigg(\langle X,Y\rangle\langle Y,U\rangle-\langle Y,Y\rangle\langle X,U\rangle\bigg)\bigg(\langle[Z,Y]_{\mathfrak{m}},V\rangle+\langle[Z,V]_{\mathfrak{m}},Y\rangle\bigg)\\
      &+2m(m+1)\langle X,Y\rangle\bigg(\langle X,Y\rangle\langle Y,V\rangle- Y,Y\rangle\langle X,V\rangle\bigg)\bigg(\langle[Z,Y]_{\mathfrak{m}},U\rangle+\langle[Z,U]_{\mathfrak{m}},Y\rangle\bigg)\Bigg].\\
     \end{split}
     \end{equation*}
     As $(G/H,\langle~,~\rangle)$ is naturally reductive, we get 
     \begin{equation*}
     g_{_{Y}}([Z,U]_{\mathfrak{m}},V)+ g_{_{Y}}([Z,V]_{\mathfrak{m}},U)+ 2C_{Y}([Z,Y]_{\mathfrak{m}},U,V)=0.
     \end{equation*}
     Hence $(G/H,F)$ is naturally reductive. 
    \end{theorem}
\section{Flag curvature of a naturally reductive homogeneous Finsler space with generalized $m$- Kropina metric}
In this section, we derive an explicit formula for the flag curvature of a naturally reductive homogeneous Finsler space with generalized $m$-Kropina metric in the sense of Deng and Hou.
\begin{theorem}
		Let $(G/H,F)$ be a naturally reductive homogeneous Finsler space with gerenalized $m$-Kropina metric $F=\dfrac{\alpha^{m+1}}{\beta^{m}}~(m\neq 0,-1)$, defined by an invariant Riemannian metric $\tilde{g}$ and an invariant vector field $\tilde{X}$ which is parallel with respect to $\tilde{g}$ such that $\tilde{X}_{H}=X$ and $\sqrt{\tilde{g}(\tilde{X},\tilde{X})}<1$. Assume that $\{U,Y\}$ is orthonormal basis for $P$ with respect to $\langle~,~\rangle$, where $(P,Y)$ is a flag on $T_{H}(G/H)$. Then the flag curvature on the flag $(P,Y)$ is given by
\begin{equation}
\begin{split}\label{a}
    K(P,Y)= \dfrac{\splitdfrac{\langle X, Y\rangle^{2m}\biggl[(m+1)\langle X,Y\rangle^{2}\bigg(\Big\langle U,\big[Y,[U,Y]_{\mathfrak{m}}\big]_{\mathfrak{m}} \Big\rangle+\Big\langle U,\big[Y,[U,Y]_{\mathfrak{h}}\big]\Big\rangle\bigg)}         {+m(2m+1)\langle X,U\rangle \bigg(\Big\langle X,\big[Y,[U,Y]_{\mathfrak{m}}\big]_{\mathfrak{m}}\Big\rangle+\Big\langle X,\big[Y,[U,Y]_{\mathfrak{h}}\Big\rangle\bigg) \bigg]}}{4(m+1)\big[m\langle X,U\rangle^{2}+\langle X,Y\rangle^{2}\big]}.
\end{split}
\end{equation}
\textit{Proof.} Since $(G/H,F)$ is naturally reductive homogeneous Finsler space, by Lemma \ref{ER}, we have 
\begin{equation}\label{gy}
R(U,Y)Y=\dfrac{1}{4}\big[Y,[U,Y]_{\mathfrak{m}}\big]_{\mathfrak{m}}+\big[Y,[U,Y]_{\mathfrak{h}}\big] ~\forall~ U,Y\in\mathfrak{m}.
\end{equation}
From equation \eqref{i} and \eqref{gy}, we get
\begin{equation*}
\begin{split}
K(P,Y)=\dfrac{\splitdfrac{\langle X, Y\rangle^{2m}\biggl[(m+1)\langle X,Y\rangle^{2}\bigg(\Big\langle U,\dfrac{1}{4}\big[Y,[U,Y]_{\mathfrak{m}}\big]_{\mathfrak{m}}+\big[Y,[U,Y]_{\mathfrak{h}}\big] \Big\rangle\bigg)}{+m(2m+1)\langle X,U\rangle\bigg(\Big\langle X,\dfrac{1}{4}\big[Y,[U,Y]_{\mathfrak{m}}\big]_{\mathfrak{m}}+\big[X,[U,Y]_{\mathfrak{h}}\big] \Big\rangle\bigg)\biggr]}}{(m+1)\big[m\langle X,U\rangle^{2}+\langle X,Y\rangle^{2}\big]}.
\end{split}
\end{equation*}
After simplification, we get equation \eqref{a}.
\end{theorem}
\begin{theorem}
Let $(G/H,F)$ be a naturally reductive homogeneous Finsler space with gerenalized $m$-Kropina metric $F=\dfrac{\alpha^{m+1}}{\beta^{m}}~(m\neq 0,-1)$, defined by an invariant Riemannian metric $\tilde{g}$ and an invariant vector field $X$ on $G$ such that $H=\{e\}$ and Chern connection of $F$ coincides with the Riemannian connection of $\tilde{g}$. Let $(P,Y)$ be a flag on $T_{e}(G)$ such that $\{U,Y\}$ be an orthonormal basis of $P$ with respect to $\langle~,~\rangle$. Then flag curvature is given by
\begin{equation}{\label{hh}}
\begin{split}
K(P,Y)=\dfrac{\langle X, Y\rangle^{2m}\bigg[(m+1)\langle X,Y\rangle^{2}\Big\langle U,\big[Y,[U,Y]\big] \Big\rangle
	+m(2m+1)\langle X,U\rangle\Big\langle X,\big[Y,[U,Y]\big]\bigg]}{4(m+1)\big[m\langle X,U\rangle^{2}+\langle X,Y\rangle^{2}\big]}.
\end{split}
\end{equation}
\textit{Proof.} Since $\langle~,~\rangle$ is bi-invariant, we have
\begin{equation*}
R(U,Y)Y=\dfrac{1}{4}\big[Y,[U,Y]\big]=R(U,Y)Y=-\dfrac{1}{4}\big[[U,Y],Y\big].
\end{equation*}
Substituting above value in equation \eqref{i}, we get
\begin{equation*}
K(P,Y)=\dfrac{\langle X, Y\rangle^{2m}\bigg[(m+1)\langle X,Y\rangle^{2}\Big\langle U,\dfrac{1}{4}\big[Y,[U,Y]\big] \Big\rangle
	+m(2m+1)\langle X,U\rangle\Big\langle X,\dfrac{1}{4}\big[Y,[U,Y]\big]\bigg]}{(m+1)\big[m\langle X,U\rangle^{2}+\langle X,Y\rangle^{2}\big]}.
\end{equation*}
After simplification, we get equation \eqref{hh}.
\end{theorem}
	\section*{Acknowledgments}
The second author is thankful to UGC for providing financial assistance in terms of JRF scholarship vide letter with Ref. No.: $961$/(CSIR-UGC NET DEC. $2018$). The third author is thankful to UGC for providing financial assistance in terms of JRF scholarship vide letter with Ref. No.: $1010$/(CSIR-UGC NET DEC. $2018$).

\end{document}